\documentclass{article}[12 pt, english]

\usepackage{amssymb}
\newcommand{\Gdag}{G^{\dagger}}

\newtheorem{Zae}{Zae}[section] 
\newtheorem{definition}[Zae]{Definition} 
\newtheorem{lemma}[Zae]{Lemma} 
\newtheorem{prop}[Zae]{Proposition}
\newtheorem{theorem}[Zae]{Theorem} 
\newtheorem{mthm}{Theorem} 
\newtheorem{coro}[Zae]{Corollary} 
\newtheorem{example}[Zae]{Example}

\newcommand{\qed}{\raisebox{-.8ex}{$\Box$}}
\newenvironment{bew}
{\noindent{\bf Proof.}}
{\hfill \qed\\}

\newcommand{\NN}{\mathbb{N}}
\newcommand{\FF}{\mathbb{F}}
\newcommand{\Sym}{~{\rm Sym}}
\newcommand{\Aut}{~{\rm Aut}}
\newcommand{\MSuz}{~{\rm MSuz}}
\newcommand{\Suz}{~{\rm Suz}}
\newcommand{\PSL}{~{\rm PSL}}

\begin{document}
\title{Finite Zassenhaus Moufang sets with root groups of even order}      
\author{Barbara Baumeister, Matthias Gr\"uninger \\ Freie Universit\"at Berlin
 \\      
Institut f\"ur Mathematik \\ Arnimallee 3 \\ 14195 Berlin \\  
E-Mail:baumeist@mi.fu-berlin.de, matgruen@zedat.fu-berlin.de}      
\maketitle      
\begin{abstract} In \cite{Suz1} Suzuki classified all Zassenhaus groups of finite odd degree.
He showed that such a group is either isomorphic to a Suzuki group or to $\PSL(2,q)$ with 
$q$ a power of $2$. In this paper we give another proof of this result using the language of Moufang 
sets. More precisely, we  show that every Zassenhaus Moufang set having root groups of finite 
even order is either special and thus isomorphic to the projective line over a finite field of 
even order or is isomorphic to a Suzuki Moufang set.

{\small {\it Keywords}: Moufang sets, Zassenhaus groups, Suzuki groups}
\end{abstract}      

\section{Introduction}

A Moufang set is a set $X$ with $|X| \geq 3$, together with a collection of subgroups $(U_x)_{x \in X}$ 
acting faithfully on $X$ (called root groups), such that each $U_x$ fixes $x$ and acts regularly on $X\setminus{\{x\}}$,
and such that $U^g_x = U_{x^g}$ for each $x \in X$ and each $g$ in $G^\dagger:= \langle U_y~|~y \in X\rangle$, 
the {\em little projective group of the Moufang set}. It is immediate from the definition that 
this group acts doubly transitively on $X$.

Moufang sets have been introduced by Tits in order to describe the absolutely simple
algebraic groups of relative rank one \cite{Ti}. The concept of a Moufang set is strongly related
to the concept of a split BN-pair of rank one. Notice that it is also closely related to the concept 
of an abstract rank one group due to Timmesfeld \cite{T}.

As usual we choose two different elements in $X$ and denote them by $\infty$ and $0$. Since $U_{\infty}$ acts 
regularly on $U:=X \setminus \{\infty\}$ it induces a unique group structure on $U$ such that $0$ is the neutral 
element and such that $U \cong U_{\infty}$ (see [5, Lemma 1.3]). The Moufang set is completely determined by 
$U$ and an arbitrary element $\tau \in \Sym(X)$ which interchanges $0$ and $\infty$ and which maps $U_{\infty}$ onto $U_0$. 
The Moufang set is also denoted by $M(U,\tau)$, see Section~\ref{notation}.  

The finite Moufang sets were classified a long time ago using different language. 
This was done by Hering, Kantor and Seitz \cite{HKS} (see also the references therein).
Their classification uses difficult and long papers such as \cite{Suz2} and \cite{GoWa}. It seems to us that the concept of a Moufang
set is the appropiate language to carry out the determination of these groups.

De Medts and Segev ([4] and [16]) gave a new proof using this language under the further condition that the Moufang
set is special -- for the definition of special see the next section. The goal of  this paper is to extend their
proof to the finite Zassenhaus Moufang sets and thereby giving a partial answer to Question 3
posed by Segev in \cite{S}. A Moufang set is {\em Zassenhaus} if 
$\Gdag$ is a  Zassenhaus group, i.e. if in $\Gdag$ there is  a non-identy element which fixes two elements in $X$,
 but only the identity  fixes three elements. 

The finite Zassenhaus Moufang sets had been determined by Feit \cite{F}, Ito \cite{I}, Higman \cite{Hi} and Suzuki \cite{Suz1} in a long proof.
There are two families of examples:
\smallskip\\
$M(\FF_q)$: The set $X$ is just the projective line ${\cal P}(q)$, $q$ a prime power,  and the little projective group is $\PSL(2,q)$ in its natural action on $X$.
\smallskip\\
$\MSuz(2^{2n+1})$: This Moufang set is the natural domain for the Suzuki group $\Suz(2^{2n+1})$ with $n \in {\mathbb N}$, see Definition \ref{Def Suz Moufang}.

In this paper we give an elementary and short proof of the classification of the finite Zassenhaus Moufang sets with
root groups of even order. The latter implies that $U$ contains an involution.
We distinguish the two cases according to  whether $U$ contains a special involution (see Definition 3.6(b)) or not.

\begin{mthm}\label{theorem1}
Let $M(U,\tau)$ be a finite Zassenhaus Moufang set such that $U$ is of even order.
If there is a special involution in $U$, then $M(U, \tau) = M(q)$ and $\Gdag \cong \PSL_2(q)$ with $q = |U|= 2^m$
for some $m$ in $\NN$.
\end{mthm}

\begin{mthm}\label{theorem2}
Let $M(U,\tau)$ be a finite Zassenhaus Moufang set such that $U$ is of even order.
If there is no special involution in $U$, then $M(U, \tau) =  \MSuz(q)$ with $q^2 = |U|$,
$q$ an odd power of $2$.
\end{mthm}

As a corollary we obtain 

\begin{coro}\label{corollary}
 Let $M(U,\tau)$ be a finite Zassenhaus Moufang set such that $U$ is of even order.
Then one of the following holds:
\begin{itemize}
 \item[(a)] $U$ is abelian, $M(U, \tau) = M(q)$ and $\Gdag \cong \PSL_2(q)$ for some even prime power $q$.
 \item[(b)] $U$ is a Suzuki $2$-group, $M(U, \tau) = \MSuz(q)$ and $\Gdag =\Suz(q)$ with $q$ an odd power of $2$.
\end{itemize}
\end{coro}

Notice that apart from  \cite{S} this paper is one of the first discussing not only special but also non-special Moufang sets.

Notice also that the distinction we make in our main theorems has in fact also been made by Suzuki
without using the language of Moufang sets.
Our proof differs heavily from Suzuki's - in particular in the case that there is a special involution in $U$.
There is also some hope that some of our arguments can be extended to the case of infinite Zassenhaus Moufang sets.

The proof of Theorem \ref{theorem1} uses only the language of Moufang sets. The strategy is to show that
$U$ is an elementary abelian $2$-group and then to quote \cite{G} or \cite{DS2}.

The proof of Theorem \ref{theorem2} is at some places a translation of the proof of Suzuki into the
language of Moufang sets. Roughly speaking the idea of the proof is first to identify the root group $U$ as a
Suzuki $2$-group $A(n, \theta)$, see \ref{Antheta}. Then to determine a partition of $U$, see \ref{partition},
 and to use this partition to determine the action of $\tau$ on $U$, see \ref{tau}.

We had some difficulties to prove that $U$ is a $p$-group. Special and non-special
 Moufang sets behave very differently. In a special Moufang set the order of every element of $U$ is
a prime number or infinity.
This is not true in an arbitrary Moufang set.
At this point we quote parts of the proof given by Feit \cite{F} and presented in \cite{HB3}. We communicated also
with Bender who believes that a  character-free proof is out of reach \cite{B}.
But we managed to show, if $U$ is a nilpotent root group of finite or infinite order, which contains for every $b\in U^\#$ a special
involution $a \in Z(U) \cap V_b$, then $U$ is an elementary abelian $2$-group and $M(U, \tau)$ is special, 
see Proposition~\ref{condspecial}.

Moreover, we refer to the classification 
of Suzuki $2$-groups in \cite{Hi} since we did not want to repeat the arguments given there, which are mainly linear algebra.
The rest of the proof is pure Moufang set theory. At some parts it is shorter and more lucid 
than the original proof of Suzuki in \cite{Suz1}. For example, in Theorem \ref{theorem2} we do not have to compute the class number of the group which turns 
out to be the Suzuki group.

The paper is organized as follows:
In Chapter 2 we introduce the notation we are using. In Chapter 3 we  present some facts about 
Moufang sets. Especially we  introduce the concept of special elements and prove some lemmata which we  use 
later. In Chapter 4 we treat Zassenhaus Moufang sets and prove Theorem \ref{theorem1}.
Then in Chapter 5 we present the definition of and some facts about (generalized) Suzuki Moufang sets, see also \cite{VM},
 and prove Theorem \ref{theorem2}. 
The proof makes use of some nice properties of the Suzuki Moufang sets, for instance that all the
involutions in a root group are conjugate. In the final chapter we discuss generalized 
Suzuki Moufang sets and see that some properties such as that all 
involutions in a root group are conjugate hold only in the 'ordinary' Suzuki Moufang sets.
\medskip\\
\noindent
{\bf Acknowledgements} We like to thank Yoav Segev for carefully reading and improving our paper as well as
Tom De Medts for his comments.

\section{Notation}\label{notation}

We first introduce some notation. 

We can obtain every Moufang set by the following procedure (see [2] or [5]). Let $(U,+)$ be a (not necessarily abelian) group. 
Set $X:=U \cup \{\infty\}$ and let $\tau$ be an element in $\Sym(X)$ which interchanges $0$ and $\infty$. 

\begin{itemize}\label{notations}
\item[(a)] For $a \in U$ let $\alpha_a$ be the map in $\Sym(X)$ defined by 
$\infty \alpha_a =\infty$ and $b\alpha_a = b+a~ \mbox{for}~ b \in U.$
\item[(b)]  Set 
$U_{\infty} =\{\alpha_a~|~ a \in U\}$ and for $a \in U$ set $U_a :=U_{\infty} ^{\tau \alpha_a}$ 
\item[(c)] For $a \in U^{\#}:= U \setminus \{0\}$ set $h_a:= \tau \alpha_a \alpha_{-a\tau^{-1}}^{\tau} \alpha_{-(-a\tau^{-1})\tau}$ 
(the {\em Hua map} corresponding to $a$). Note that $h_a$ fixes $\infty$ and $0$.
\end{itemize}

Then $(X, (U_x)_{x \in X})$ is a Moufang set iff the restriction of each Hua map to $U$ is contained
in $\Aut(U)$, see [5, 3.2].  We set $$M(U,\tau):= (X, (U_x)_{x \in X}).$$ As usual for $x \in X$ we denote by $\Gdag_x$
the stabilizer of $x$ in $\Gdag$.

A Moufang set $(X,(U_x)_{x \in X})$ is called {\em proper} if $U_x \ne \Gdag_x$ for $x \in X$, or equivalently, if 
$\Gdag$ is not sharply $2$-transitive on $X$.

We recall some more standard notation which can for instance be found in [2].

\begin{itemize}\label{notations2}

\item[(d)]  For $a \in U^\#:= U\setminus{\{0\}}$, let $\mu_a =\alpha_{(-a)\tau^{-1}}^{\tau} \alpha_a \alpha_{-a\tau^{-1}}^{\tau} $
 be the unique element in $U_0 \alpha_a U_0$ with $\infty \mu_a = 0$ and $0 \mu_a = \infty$ (see \cite{DS1}, 4.1.1). Note that 
$M(U,\tau) =M(U,\mu_a)$. 

\item[(e)] Set $H:= \langle \mu_a \mu_b~|~ a,b  \in U^{\#} \rangle$, the {\em Hua subgroup of $M(U,\tau)$}. 
By definition $H \leq \Gdag_{0,\infty}$ and thus $H \leq N_{\Gdag} (U_{\infty})$. Therefore $H$ acts as automorphism group 
on $U$. 
\end{itemize}

As in [2, Proposition 2.10] we write:

\begin{enumerate}
\item[(f)] For $a \in U^\#$ set $\sim a := (-a\tau^{-1})\tau$. One easily sees $\sim (\sim a) = a$. 
In particular, $\sim$ is bijective.
\end{enumerate}

\section{Preliminary observations}

In the following section, $M(U,\tau)$ is an arbitrary, not necessarily finite Moufang set. 
We first recall some properties of the $\mu$-functions and $\sim$ which we will frequently use. 
All of them are already known; we only add a proof if we didn't find it explicitely stated anywhere.

\begin{lemma}\label{mu and sim}
Let $a \in U^{\#}$ and $h \in H$. Then
\begin{enumerate}
\item $\mu_a^{-1} = \mu_{-a}$
\item If $M(U,\tau) =M(U,\tau^{-1})$, then $\mu_{a\tau} =\mu_{-a}^{\tau}$. 
\item $\mu_{a\mu_b} = \mu_{-a} ^{\mu_b}$.
\item $H=\Gdag_{0,\infty}$.
\item $\sim a = -(-a)\mu_a$. Especially, this implies that the element $\sim a$ does not depend on the choice of $\tau$.
\item $\mu_{-a} = \alpha_{-\sim a} \mu_{-a} \alpha_a \mu_{-a} \alpha_{\sim -a}$.
\item $\mu_{ah} = \mu_a ^h$.
\item $\sim (ah) = (\sim a) h$.
\item $\mu_{\sim a} = \mu_{-a}$ and $ \mu_a = \mu_{-\sim a} = \mu_{\sim - a}$.
\item $\sim (a\tau) = (-a)\tau$ and $-(a\tau^{-1}) = (\sim a)\tau^{-1}$.
\item $a\mu_a=\sim-\sim a$ and $a\mu_{-a} =-\sim -a$. Especially $-\sim -a = \sim -\sim a$ 
if $\mu_a$ is an involution.
\item If $b \in U^{\#}$ with $a\ne b$ and $\mu_a =\mu_b$, then
$\mu_{a\mu_a^{-1} -b\mu_a^{-1}} = \mu_{a-b} $ and $\mu_{-a+b} =\mu_{-\sim a+ \sim b}$.
\end{enumerate}
\end{lemma}
\begin{bew}
Parts (a)-(c) and (e)-(g) follow from proposition 4.3.1 of \cite{DS1}. Part (d) follows from 4.2.2 of 
[3]. (k) and (l) can be found in 2.3 and 2.5 of [15]. 

(h) Since $h$ is an automorphism of $U$, we get with (e) and (g) that 
$$\sim (ah) = -(-ah)\mu_{ah} = -(-ah) h^{-1} \mu_a h=
-(-a)\mu_a h = (\sim a) h.$$

(i) We have by  parts (e) and (b) that 
$$\mu_{\sim a} = \mu_{-(-a)\mu_a} = (\mu_{(-a)\mu_a})^{-1} =
((\mu_a)^{\mu_a})^{-1}= \mu_a^{-1} = \mu_{-a},$$
the first part of the assertion. This yields $\mu_{-\sim a} = \mu_{\sim a}^{-1} = \mu_{-a}^{-1} =\mu_a$ as well as 
$\mu_{\sim -a} = \mu_{--a} = \mu_a$.

(j) By definition $$\sim (a\tau) =(- ((a\tau)\tau^{-1}))\tau = (-a)\tau$$ and 
$$-(a\tau^{-1})=((-(a\tau^{-1}))\tau\tau^{-1} = (\sim a)\tau^{-1}.$$

\end{bew}

We will repeatedly use the following two fundamental equations (see \cite{DS1}, 6.1.1):
\begin{itemize}
\item[(3A)] If $a,b \in U^{\#}$ with $a \ne b$, then the element $c:=(a\tau^{-1} -b\tau^{-1})\tau $ does not depend on $\tau$.
More precisely, $c=(a-b)\mu_b + \sim b$.
\item[(3B)] $\mu_c = \mu_{-b} \mu_{b-a} \mu_a$.
\end{itemize}

\subsection{Some properties of involutions in a root group}

The $\mu$-maps are very important in the theory of Moufang sets. Therefore it is helpful to know the following.
Part (a) is already included in \cite{DS2} 7.3.1 (10) in the case of special Moufang sets.

\begin{lemma}\label{a inv}
If $a$ is an involution in $U^{\#}$, then 
\begin{enumerate} 
\item $\mu_a^{\alpha_{-\sim a}} = \alpha_a^{\mu_a}$. Especially $\mu_a$ is an involution which is conjugate to $\alpha_a$.
\item $\sim a$ is the unique fixed point of $\mu_a$. 
\item  $\sim - \sim a =a\mu_a=-\sim a$. 
\end{enumerate}
\end{lemma}
\begin{bew} 
(a)  By \ref{mu and sim} (f) and as $a = -a$ we have  $\mu_a = \alpha_{-\sim a} \mu_a \alpha_a \mu_a \alpha_{\sim a}=
\alpha_a^{\mu_a \alpha_{\sim a}}$.

(b) Since $\infty$ is the unique fixed point of $\alpha_a$, we have that $\infty{\mu_a \alpha_{\sim a}} =
0{\alpha_{\sim a}} = \sim a$ is 
the unique fixed point of $\mu_a = \alpha_a^{\mu_a \alpha_{\sim a}}$.

(c) The assertion follows with \ref{mu and sim} (k) as $a = -a$.
\end{bew}

In general, different elements in $U^{\#}$ can induce the same $\mu$-map. As in [15, equation (1.1)] we denote 
$V_a:=\{ b \in U^{\#}~|~\mu_b =\mu_a\}$ 
for $a \in U^{\#}$. Notice that $-V_a = V_{-a}$ by \ref{mu and sim} (a) and that
$-\sim a,\sim -a \in V_a$ by \ref{mu and sim} (i). We next show that each of these sets contains at most one 
involution (see also \cite{DS2} 7.3.1 (6) in the case of special Moufang sets).

\begin{lemma}\label{mua = mub}
If $a,b \in U$ are involutions with $\mu_a = \mu_b$, then
$a=b$.
\end{lemma}
\begin{bew}  If $\mu_a =\mu_b$, then by Lemma \ref{a inv} (b) we get $\sim a = \sim b$ and thus
$a=b$. 
\end{bew}

If $U$ is a finite group of even order, then the involutions in $U$ behave nicely, as we show next.

\begin{lemma}\label{order mua mub}
Let $M(U,\tau)$ be a Moufang set such that $U$ has finite 
even order. Then $\mu_a \mu_b$ has odd order for all involutions $a,b \in U$.
Hence all involutions in $U$ are $H$-conjugate.
 \end{lemma}
\begin{bew} 
 We prove the first statement by induction on $|U|$. Suppose 
$a, b \in U$ are involutions such that $\mu_a \mu_b$ has even order $2n$.
Set $t:=(\mu_a \mu_b) ^n$. Then $t \in H$ and $t$ centralizes $\mu_a$ 
and $\mu_b$. It follows that $t$ centralizes $ a$ and $b$ as well. Hence 
$a,b \in V:=C_U(t)$ which is a root subgroup of $U$ (\cite{DS1}, 6.2.3). By 6.2.2 of the same paper one can choose $\tau$ in  
such a way that $\tau$ stabilizes $V^{\#}$ and that if $\rho$ is the restriction of $\tau$ to $V \cup \{\infty\}$, then 
$M(V,\rho)$ is a Moufang set. Since $t \ne 1$, $V$ is 
a proper subgroup of $U$ and hence we can apply the induction hypothesis. Hence there is an 
odd number $k$ such that $(\mu_a \mu_b) ^k $ centralizes $V$. Set $l ={{k+1} \over 2}$ and $h =
(\mu_a \mu_b)^l$. Thus we get $$h^2 |V = (\mu_a \mu_b)^{2l}|V =(\mu_a \mu_b)^{k+1}|V = \mu_a \mu_b|V$$ and 
$$\mu_{ah}|V^{\#}=h^{-1} \mu_a h|V^{\#} =
\mu_a h^2|V^{\#} = \mu_a \mu_a \mu_b |V^{\#} = \mu_b|V^{\#}.$$ 
Using Lemma \ref{mua = mub} for $M(V,\rho)$, this implies $ah = b$. But then 
we get $\mu_a ^h = \mu_{ah} = \mu_b$ and thus $h^{-2} \mu_a \mu_b = \mu_a^h \mu_b = 1$. 
Therefore $(\mu_a \mu_b)^{2l} = h^2 =\mu_a \mu_b$ and thus $(\mu_a \mu_b)^{2l-1} =1$.
This contradicts our assumption that $\mu_a \mu_b$ has even order. 
Hence we have proved the first statement. 

The second statement follows 
immediately since we have shown that $\mu_a $ and $\mu_b$ are $H$-conjugate for 
all involutions $a,b \in U$, which together with \ref{mua = mub} implies that $a$ and $b$ are $H$-conjugate as 
well.
 \end{bew}

We remark here that in the infinite case it is possible that there is more than one $H$-orbit of involutions in a root group. This happens for example in ${\mathbb M}(K)$ if $char K=2$, $K$ not perfect, or in $\MSuz(K,L,\theta)$ for $\theta$ not surjective (see Section 6).

The following lemma is useful for determining the $\mu$-maps. We will apply it in \ref{orbits2}. 
The given proof is due to Segev who  simplified our original proof.

\begin{lemma}\label{orbit of subgroup of H}
Let $M(U,\tau)$ be a Moufang set with Hua subgroup $H$ and 
let $V$ be a subgroup of $U$. Suppose that

 \begin{enumerate}
\item there is an abelian subgroup $K$ of $H$ such that all elements in
 $V^{\#}$ are $K$-conjugate,
\item $h^{\mu_a} = h^{-1}$ for all $a \in V^{\#}, h \in K$.
\end{enumerate}

 Then $\mu_{-a+b} =\mu_{ -\sim b +\sim a  }$ for all 
$a,b \in V^{\#}$ with $a \ne b$. 
\end{lemma}

\begin{bew} 
There exist elements $g,h \in K$ with 
$ah = b, ag = b-a$.
We compute 
$$(a\mu_{-a} -b\mu_{-b}) \mu_a =(a\mu_{-a} -ah\mu_{-ah})\mu_a = (a\mu_{-a} -ah \mu_{-a}^h)\mu_a=$$
$$(a\mu_{-a} -ah h^{-1} \mu_{-a} h) \mu_a = (a\mu_{-a}- ah^{-1} \mu_{-a} ) \mu_a.$$
Further by (3B), 
$$\mu_{(a\mu_{-a} -ah^{-1}\mu_{-a})\mu_a} = \mu_{-ah^{-1}} \mu_{ah^{-1} -a} \mu_a =h \mu_{-a} h^{-1} \mu_{-agh^{-1}} \mu_a=$$
$$h^2 \mu_{-a} hg^{-1} \mu_{-a} gh^{-1} \mu_a=h^2 \mu_{-a} h^2 g^{-2} \mu_{-a} \mu_a = g^2 \mu_{-a}.$$
It follows that 
$$\mu_{(a\mu_{-a} -b\mu_{-b})\mu_a} =g^2 \mu_{-a}.$$
In particular $\mu_{(a\mu_{-a} -b\mu_{-b})\mu_a}$ inverts all elements of $K$ and so 
$$\mu_{-a} = g^{-2} \mu_{(a\mu_{-a} -b\mu_{-b})\mu_a} = g^{-1}  \mu_{(a\mu_{-a} -b\mu_{-b})\mu_a} g =
g^{-1} \mu_{b\mu_{-b} -a\mu_{-a}}^{\mu_a} g =$$
$$ \mu_{b\mu_{-b} -a\mu_{-a}}^{\mu_a g} =\mu_{b\mu_{-b} -a \mu_{-a}}^{g\mu_{-a}}.$$
It follows that 
$$\mu_{b\mu_{-b} -a\mu_{-a}} = \mu_{-a}^{\mu_a g} = \mu_{-a}^g = \mu_{-ag} = \mu_{a-b}.$$
Now $b\mu_{-b} -a\mu_{-a} = -\sim -b +\sim -a$ by 3.1 (k). If we replace $a$ by $-a$ and $b$ by $-b$, we finally get the 
desired formula. \end{bew}

\subsection{Special elements}

The concept of special Moufang sets is well established. 
In this section we introduce the concept of a special element in a root group of a Moufang set, which in fact already appeared
in  \cite{T}, $\S 2$. 

\begin{definition}\label{special}
\begin{enumerate}
\item A Moufang set $M(U,\tau)$ is called special if $(-a)\tau =-(a\tau)$ for all $a \in U^{\#}$.
\item 
An element $a \in U^{\#}$ is called special if
$(-a){\tau^{-1}} =-(a{\tau^{-1}})$. 
\end{enumerate}
\end{definition}

It is easy to see that $M(U,\tau)$ is special iff all the elements in $U^{\#}$ are special:
If $M(U,\tau)$ is special and $a \in U^{\#}$, then 
$$(-a)\tau^{-1} = (-(a\tau^{-1} \tau))\tau^{-1} =
(-(a\tau^{-1}))\tau \tau^{-1} =-(a\tau^{-1}).$$ 
If all elements in $U^{\#}$ are special, then for $a \in U^{\#}$ we get
$$(-a)\tau = (-((a\tau)\tau^{-1}))\tau = (-(a\tau))\tau^{-1}\tau =-(a\tau).$$ 

It might surprise the reader that we require $(-a){\tau^{-1}} =-(a{\tau^{-1}})$ and not 
$(-a){\tau} =-(a{\tau})$ for the definition of a special element. But since 
$M(U,\tau) =M(U,\phi \tau)$ for all $\phi \in \Aut (U)$, it may happen
that there are $\tau, \rho \in \Sym(X)$ and $a \in U^{\#}$ with $M(U,\tau) = M(U,\rho)$ and $ (-a)\tau =-a\tau$ but 
$(-a)\rho \ne -a\rho$.  An example for this situation is a Ree-Tits Moufang set 
where there are special elements whose conjugates are not always special. 
Whereas in the next lemma we show that $(-a)\tau^{-1} =
-(a\tau^{-1})$ iff $(-a)\rho^{-1} = -(a\rho^{-1})$ and so our definition of special 
does not depend on the choice of $\tau$.
Moreover, if $M(U,\tau) \ne M(U,\tau^{-1})$, then $(-a){\tau} =-(a{\tau})$ is in general not equivalent to 
one of the statements in the following lemma. 

\begin{lemma}\label{equivalence special} 
For $a \in U^{\#}$ the following statements are equivalent.
\begin{enumerate}
\item $a$ is special.
\item $\sim a= -a$.
\item $(-a) \mu_a = a$.
\item $a \mu_{-a} = -a$.
\item If $M(U,\tau) =M(U,\rho)$, then $(-a){\rho^{-1}} =-(a{\rho^{-1}})$.
\item There is an element $\rho \in \Sym(X)$ with $M(U,\rho) =M(U,\tau)$  
such that $-(a{\rho^{-1}}) =(-a){\rho^{-1}}$.
\item $(-a)\mu_a = -(a\mu_a)$.
\end{enumerate}
\end{lemma}

 \begin{bew}
(a) implies (b): By definition, we have $\sim a = 
(-(a\tau^{-1}))\tau = ((-a)\tau^{-1})\tau = -a$.

(b ) and (c) are equivalent: We have by \ref{mu and sim} (e) 
$\sim a = -(-a)\mu _a$ and thus $\sim a= -a$ iff ($-a)\mu_a = a$.

(c) and (d) are equivalent: This is clear since $\mu_a^{-1} = \mu_{-a}$.

(b) implies (e): By the definition of $\sim a$ we have $-a = \sim a = (-a \rho^{-1})\rho$ and hence 
$(-a)\rho^{-1} = -a\rho^{-1}$.

(e) implies (f): This is trivial.

(f) implies (g): Again by 3.5 and 4.4.1 (1) in \cite{DS1} $\rho \mu_a$ 
induces an automorphism of $U$. 
Therefore, we have $(-a)\mu_a = (-a)\rho^{-1} \rho \mu_a =((-a) \rho^{-1}) \rho
\mu_a = (-a\rho^{-1}) \rho\mu_a = -(a\rho^{-1}) \rho \mu_a = -(a\mu_a)$.

(g) implies (a): By 4.3.1 (1) in \cite{DS1}, $\mu_a =\tau^{-1} h_a$ 
where $h_a$ is the Hua map associated to $a$ and hence induces an automorphism 
on $U$. Thus $-(a\tau^{-1}) = -a (\mu_a h_a^{-1}) = (-(a\mu_a)) h_a^{-1} = 
(-a)\mu_a h_a^{-1} = (-a) \tau^{-1}$. 
\end{bew}

Here are some more properties of special elements.

\begin{lemma}\label{prop of special}
\begin{enumerate}
\item An element $a \in U^{\#}$ is special iff $-a$ is
special. 
\item If $a \in U$ is an involution, then $a$ is special iff 
$a\tau^{-1}$ is again an involution.
\item If $a \in Z(U)^{\#}$ is special, then $a \rho^{-1} \in Z(U)$ for all 
$\rho \in \Sym(X)$ with $M(U,\tau) =M(U,\rho)$.
\end{enumerate}
\end{lemma}

\begin{bew} The first statement follows from the fact that (c) and (d)  of \ref{equivalence special} 
are equi\-valent, the second is true by definition. If $a \in Z(U)^{\#}$ is 
special and $\rho$ as above, then $\mu_a \rho^{-1}$ induces an automorphism of $U$.
Since $a =(-a)\mu_a$ we get $a\rho^{-1} =(-a)\mu_a \rho^{-1} \in Z(U)$. 
\end{bew}

\begin{lemma} \begin{enumerate}
\item
If $a \in U^{\#}$ is special, then $a\mu_a =-a = a\mu_{-a}$. 
\item If $a$ has order $4$, then $a$ is not special. 
\end{enumerate}
\end{lemma}
\begin{bew} 
\begin{enumerate} 
\item By \ref{equivalence special} (c) and (g), we have $-(a\mu_a) = (-a)\mu_a = a$ and thus $a\mu_a =-a$.
The second equation holds by \ref{equivalence special} (d).
\item Suppose $a$ has order $4$ and is special. Then with (3A) and part (a)
$$(a \cdot 2) \mu_{-a}+ a =(a-(-a))\mu_{-a} +\sim -a = $$
$$(a\mu_a -(-a)\mu_a)\mu_{-a} = (-a-a)\mu_{-a} =(a\cdot 2)\mu_{-a}.$$
But this implies $a=0$, a contradiction. 
\end{enumerate} 
\end{bew}

It is not clear whether $a\mu_a =-a$ implies that $a$ is special. Note also that the Moufang sets of  Ree-Tits type 
contain special elements of order $9$. 

To be special is an $H$-invariant property:

\begin{lemma}\label{ah special}
If $a \in U$ is special and $h \in H$, then 
$ah$ is special.
\end{lemma}

\begin{bew} If $a$ is special, then by \ref{mu and sim} (g) we have $\sim (ah ) = (\sim a) h = 
(-a) h = -ah$, hence $ah$ is special. 
\end{bew}

\begin{lemma}\label{cond special}
An element $a \in U^{\#}$ is special iff there is an element
$b \in U_0$ such that $\mu_a = b \alpha_a b$. This element is $b=\alpha_a^{\mu_a}$.
\end{lemma}

\begin{bew} We have 
$\mu_a = b^{\prime} \alpha_a b^{\prime \prime}$ with 
$b^{\prime} = \alpha_{(-a)\tau^{-1}}^{\tau} $ and $b^{\prime \prime} =
\alpha_{-a\tau^{-1}}^{\tau}$. Thus $a$ is special iff these two elements are equal.
In this case we get $b^{\prime} =b^{\prime \prime} = \alpha_a^{\mu_a}$ for 
$\tau =\mu_a$.
\end{bew}

Notice that, in fact, \ref{cond special} is (2.2) of \cite{T}.

The following lemma collects some useful information about $V_a$ for a special central element $a$.
We will need only part (a) - (d), but the other parts yielding the fact that $a$ and $-a$ are the
only special elements in $V_a$ are interesting too. This lemma should be compared to Lemma 3.1 in 
\cite{S} from which it was inspired. 

\begin{lemma}\label{a in center special} 
Suppose that $M(U,\tau) = M(U,\tau^{-1})$, that $a \in Z(U)^{\#}$ is  special, that $\mu_a =\mu_{-a} =
\mu_a^{-1}$ 
and that $b \in V_a \setminus \{a,-a\}$. Then the following hold:
\begin{enumerate}
\item $-(b-a)\mu_a+(a-b)\mu_a= \sim -b +a -\sim b$.
\item $-(a-b)\mu_a+(b-a)\mu_a= b+a\cdot 2$.
\item $-a\cdot 3 =\sim -b -\sim b +b= -\sim b + b +\sim - b= b+\sim- b -\sim b$.
\item $-((-a)\tau^{-1} -(-b)\tau^{-1})\tau+(a\tau^{-1} -b\tau^{-1} )\tau  = a$
and $(a\tau^{-1} -b\tau^{-1})\tau -((-a)\tau^{-1} -(-b)\tau^{-1})\tau  = a$ .
\item $(a-b)\tau -(-a-\sim b)\tau = a\tau$.
\item $- (-b\tau^{-1} -a\tau^{-1})\tau + ((-b)\tau^{-1}-a\tau^{-1})\tau =
 -\sim b -a$.
\item $-(-a-b)\tau+(\sim b -a)\tau =-\sim(b\tau) - a\tau$.
\item $-(-a -b)\mu_a+(\sim b -a) \mu_a= \sim b +a$.
\item $a$ and $-a$ are the only  special elements in $V_a$.
\end{enumerate}
\end{lemma} 

\begin{bew}
(a) By (3A) we have 
$$ (a-b) \mu_a+ \sim b = 
(a\tau^{-1} -b\tau^{-1} )\tau = (-b\tau^{-1}-( -a)\tau^{-1}) \tau =$$
$$((\sim b) \tau^{-1} -(-a)\tau^{-1}) \tau = (\sim b +a) \mu_a+a,$$
hence
$$(*) \ -(\sim b+ a) \mu_a+ (a-b) \mu_a= a - \sim b.$$ 
Furthermore $$(**) \ (\sim b +a )\mu_a= (-(-b)\mu_a+(-a)\mu_a)\mu_a= 
((-a)\mu_a- (-b)\mu_a)\mu_a=$$ 
$$(b-a) \mu_a+ \sim -b
=-(-\sim - b -(b-a) \mu_a).$$
If we insert $(**)$ in $(*)$, the claim follows. 

(b) This is similar to (a): We have again by (3A)
$$(b-a)\mu_a=(b\tau^{-1} -a\tau^{-1}) \tau +a = ((-a)\tau^{-1} --b\tau^{-1})
\tau +a =$$
$$((-a)\tau^{-1} -(\sim b)\tau^{-1})\tau +a = (-a-\sim b)\mu_a+ b+a ,$$
hence
$$-(-a-\sim b)\mu_a+(b-a)\mu_a= b+a.$$
Furthermore 
$$(-a -\sim b)\mu_a= ((-b)\mu_a-(-a)\mu_a)\mu_a= (-b +a)\mu_a+a =(a-b)\mu_a+ 
a.$$
Thus the claim follows.

(c) By (a) and (b), 
$$\sim - b + a -\sim b =  - a \cdot 2-b.$$
Hence $$\sim -b -\sim b +b  =-a \cdot 3.$$
We get the other equations if we conjugate by $\sim - b$ and $- \sim b$.

(d) We have with (3A) 
$$(a-b) \mu_a=(a\tau^{-1} -b\tau^{-1} ) \tau -\sim b$$
and
$$(b-a)\mu_a =  (-a-(-b))\mu_a = ((-a)\tau^{-1} -(-b)\tau^{-1})\tau -\sim - b.$$
By part (a), we get 
$$\sim -b + a -\sim b = -(-a-(-b)) \mu_a+ (a-b)\mu_a= $$
$$\sim -b -((-a)\tau^{-1} -(-b)\tau^{-1}) \tau +(a\tau^{-1} -b\tau^{-1})\tau 
- \sim b.$$
Therefore, we get the first equation. The second follows since $a$ is central 
and we can conjugate with $-(a\tau^{-1} -b\tau^{-1})\tau$.

(e)
This follows by (d) and by replacing $b$ with $b\tau$ and $a$ with 
$a\tau$.

(f)
 We have 
$$(-a\tau^{-1}-b\tau^{-1})\tau = (-a-b)\mu_a+ \sim b=-(-\sim b-(-a-b)\mu_a) $$
and 
$$((-b)\tau^{-1} -a\tau^{-1})\tau =(-a-b)\mu_a- a.$$
Subtracting yields the result.

(g) 
This follows from (f) by taking $a\tau$ instead $a$ and $b\tau$ of $b$.

(h)  We have $a\mu_a =-a$ and $b\mu_a =b\mu_b =\sim-\sim  b$ by \ref{mu and sim} (k).
Thus $-\sim (b\mu_a) = -\sim (\sim -\sim b)  = \sim b$ and the claim follows from 
(g) with $\tau =\mu_a$.

(i) If $b$ is special, then $\sim b = -b$ and hence $a=b$ by (h), a contradiction.
\end{bew}

\begin{prop}\label{mua = mu(a+b)}
If $a$ and $b$ are as in \ref{a in center special}, then $\mu_a= \mu_{a-b}
 \mu_{a\cdot 5+b} \mu_{a-b}$. Especially, if $a$ is an involution, then 
$\mu_a=\mu_{a+b}$. 
\end{prop}
\begin{bew} Set $x:= a\mu_a-b\mu_a$ and $y:= -a\mu_a-(-b)\mu_a=(-a)\mu_a -(-b)\mu_a$.
 Then $\mu_{a-b}
 =\mu_x$ and 
$\mu_y = \mu_{-a+b} =\mu_{b-a}= \mu_{a-b}^{-1}$ by \ref{mu and sim} (l). Furthermore, 
(3B) tells us 
$$\mu_c = \mu_{-y} \mu_{y-x} \mu_x$$ with
$c=(x\mu_a-y\mu_a)\mu_a$. 
We have $c= a\mu_a=-a$ by \ref{a in center special} (d) with $\tau =\mu_a$, and 
$$y-x = -a\mu_a-(-b)\mu_a-a\mu_a+b\mu_a= a + \sim b +a -\sim - b =$$
$$a\cdot 2 -   (\sim -b -\sim b) = a\cdot 2 -(-a \cdot 3 -b) = a \cdot 5 +b$$ 
by \ref{a in center special} (c).
Hence $\mu_a=\mu_{-a}= \mu_{a-b} \mu_{a\cdot 5 +b} \mu_{a-b}$. If $a$ is an involution, 
then $a-b = -(a+b)$ and hence $\mu_a=\mu_{a-b} \mu_{a+b}  \mu_{a+b}^{-1} = \mu_{a-b}$.
As $\mu_a$ is an involution by assumption, it follows that $\mu_a =  \mu_{a-b} =\mu_{a+b} $
\end{bew}

The following lemma is to some extent the converse of Proposition \ref{mua = mu(a+b)}.

\begin{lemma}\label{a in center special II} 
If $a\in Z(U)^{\#}$ is a special involution and $x \in U^{\#}$ 
with 
$\mu_{x+a} = \mu_x=\mu_{-x}$, then $\mu_x=\mu_a$.
\end{lemma}
\begin{bew} We have $((x+a)\tau^{-1} -x\tau^{-1})\tau = a\mu_x +\sim x$ and 
thus (3B) implies 
$$ \mu_{a\mu_x}= \mu_x \mu_a \mu_x = 
\mu_x \mu_a \mu_{x+a} = \mu_{a\mu_x + \sim x}.$$ 
Now $a\mu_x = a\mu_a \mu_x$ is again a special involution in $Z(U)$. 
Thus by \ref{mua = mu(a+b)} with $b = a\mu_x + \sim x$
 $$\mu_{a\mu_x } = \mu_{a\mu_x + \sim x +a\mu_x} =\mu_{\sim x} =\mu_x.$$
Hence $\mu_x \mu_a \mu_x = \mu_x$ and therefore $\mu_a =\mu_x.$
\end{bew}

Proposition~\ref{condspecial} ahead is the key for the  proof of Theorem~\ref{theorem1},
but first we need further knowledge on nilpotent groups.

\begin{lemma}\label{nilpotent} Let $G$ be a nilpotent group. Then the following holds:
\begin{enumerate}
\item If there is an element $x \in G$ of order $p$, $p$ a prime number, then there is an element
$y \in Z(G)$ of order $p$.
\item If there is an element $x\in G$ of infinite order and if there is a natural number $e \geq 1$ such that
all the elements in $G$ of finite order have order at most $e$, then there is an element $y$ in $Z(G)$ which has infinite order.
\end{enumerate}
\end{lemma}
\begin{bew}
(a) 
We prove this by induction on the nilpotency class of $G$. Let $y \in G$ have order $p$. 
If $y \not\in Z(G)$, then $G/Z(G)$ has an element of order $p$. 
Since the claim holds for $G/Z(G)$, there is $xZ(G) \in Z(G/Z(G))$ with $o(xZ(G))=p$. 
This means that $x$ is not in $Z(G)$, but $x^p \in Z(G)$.
Hence there is an element $g  \in G$ with $[x,g] \ne 1$. Moreover $[x,g] \in Z(G)$
and $[x,g]^p = [x^p,g] = 1$. Thus $y:= [x,g] $ is the desired element of order $p$ in 
$Z(G)$.

(b) Suppose that there is an element $x \in G$ such that the element $xZ(G)$ in $G/Z(G)$ 
has finite order $f$ with $f > e$. Then $x$ has infinite order and so has 
$x^f$ which is an element in $Z(G)$. Thus we can assume that all elements in $G/Z(G)$ with finite 
order have order at most $e$. If all elements in $G/Z(G)$ have finite order, then $Z(G)$ must 
contain an element of infinite order since $G$ does. Therefore we can apply induction and assume that there is an element 
$xZ(G) \in Z(G/Z(G))$ with infinite order. Again for all $g \in G$ the map $[.,g]:\langle Z(G),x \rangle \to 
Z(G): h \mapsto [h,g]$ is a homomorphism. If there is a $g \in G$ such that the image of $[.,g]$ has infinite 
order, then we are done. If not, then the element $x^f$ with $f:=e!$ is in the kernel of $[.,g]$ for all 
$g \in G$ and thus $x^f \in Z(G)$.
\end{bew}

The second condition in part (b) is essential. One can readily construct an infinite nilpotent group such that all 
elements of infinite order are not contained in the center. 

The following proposition, which holds for finite as well as for infinite Moufang sets, is related to Theorem~C of 
\cite{DST}. But the proofs are completly different as we do not assume $M(U,\tau)$ to be special.

\begin{prop}\label{condspecial}
 Suppose that $M(U,\tau)$ is a proper Moufang set such that
\begin{enumerate}
\item[(a)] $U$ is nilpotent.
\item[(b)] For every $b \in U^{\#}$ there is a special involution $a \in Z(U) \cap V_b$.
\end{enumerate}
Then $U$ is an elementary abelian $2$-group and $M(U,\tau)$ is special.
\end{prop}

\begin{bew}
We first note that since $U$ is proper there are at least two involutions in $U$.
By \ref{mua = mub} every involution in $U$ is central and special.
Suppose that $U$ is not of exponent $2$. Since $U$ is nilpotent, this means by \ref{nilpotent} that either $U$ contains an 
element of order $4$ or an element in $Z(U)$ which has odd or infinite order.

Suppose $U$ contains an element $b$ of order $4$. 
By assumption and by \ref{mua = mub} there is a unique involution $a$ in $V_b$. Then $\mu_b =\mu_b^{-1}= 
\mu_{-b} = \mu_{b + b\cdot 2}$, hence $a = b\cdot 2$ by \ref{a in center special II}. If $t$ is an involution 
distinct from $a$, then again $(b+t)\cdot 2 = b\cdot 2 = a$, as $t$ is in $Z(U)$. Hence $\mu_{b+t} = \mu_a 
=\mu_b$ and so $\mu_a=\mu_b=\mu_t$ by \ref{a in center special II} and $a=t$ by \ref{mua = mub}, a contradiction. 
Thus $U$ does not contain elements of order $4$.

We may thus assume  that there is an element $b \in Z(U)$ whose order is not a power of $2$.
Let $a$ be the unique involution in $V_b$. By \ref{a in center special} (c), 
$a = \sim -b -\sim b +b$, 
hence $a-b =\sim-b -\sim b$. Since $a-b \in Z(U)$, we get 
$a-b = -\sim b + \sim -b$. By \ref{mu and sim} (k) and as 
$$(*) \hspace{1cm} \mu_b = \mu_a = \mu_a^{-1} = \mu_b^{-1} = \mu_{-b}$$
this implies $a-b = (-b)\mu_a - b\mu_a$. 

By \ref{mu and sim} (l) $\mu_{(-b)\mu_{-b}^{-1} - b\mu_b^{-1}} = \mu_{-b -b}$. Thus by applying again $(*)$
 we see that $\mu_{a-b} = \mu_{(-b)\mu_a -b\mu_a} =\mu_{-b -b} = \mu_{-b \cdot 2} =\mu_{b\cdot 2}^{-1}$.

As $[a,b] = 1$, we have $-(a+b) = a-b$ and therefore $\mu_{a+b} =\mu_{a-b}^{-1} = \mu_{b \cdot 2}$.
With \ref{mua = mu(a+b)} we finally get $\mu_b=\mu_a =\mu_{a+b} = \mu_{b\cdot 2}$.

But again, if $t \in U$ is an involution distinct from $a$, then we get by replacing $b$ by $b+t$ that 
$\mu_{b+t} = \mu_{(b+t)\cdot 2} =\mu_{b\cdot 2 } = \mu_b$, hence $\mu_a =\mu_b=\mu_{b+t} =\mu_t$ by \ref{a in center special II}. 
By \ref{mua = mub}, this is a contradiction.
Thus $U$ is of exponent $2$. Since all elements of $U^{\#}$ are special, $U$ is special.
\end{bew}

\section{Zassenhaus Moufang sets}

In this section we discuss general facts about Zassenhaus Moufang sets and prove Theorem~\ref{theorem1}.

\subsection{Elementary facts about Zassenhaus Moufang sets}

\begin{definition} A proper Moufang set $M(U,\tau)$ is called a Zassenhaus 
Moufang set if $\Gdag_{0,\infty, a} =1$ for all $a \in U^{\#}$.
\end{definition} 

It is easily seen that $M(U,\tau)$ is a Zassenhaus Moufang set exactly if 
$C_U(h)=1$ for all $h \in H^\#$.
 
From now on, we assume that $M(U,\tau)$ is a Zassenhaus Moufang set and
that the order of $U$ is finite.

\begin{prop}\label{properties U} 
\begin{enumerate} 
\item The root group $U$ is nilpotent. 
\item If $U$ is abelian, then $M(U,\tau) \cong M({\mathbb F}_q)$ 
for $q=|U|$ and hence $\Gdag \cong \PSL_2(q)$. 
\item $\Gdag$ is simple.
\end{enumerate} 
\end{prop}

\begin{bew}
(a)
Since $M(U,\tau)$ is proper, $H \ne 1$. Thus $UH$ is a Frobenius 
group with Frobenius kernel $U$. By Thompson's theorem (\cite{Hu}, V.8.7), $U$ is nilpotent.

(b) By the main theorem of \cite{S}, $M(U,\tau)$ is special. Thus the 
claim follows with \cite{DS2} and \cite{S1}. 

(c) Suppose $1 < M$ is a normal subgroup of $\Gdag$. Since $\Gdag$ acts 
primitively on $X$, $M$ is transitive on $X$ and thus on the set of root groups.
By definition $\Gdag$ is generated by the root groups and hence we have 
$\Gdag =MU_{\infty}$ and $\Gdag/M \cong U_{\infty}/(M \cap
U_{\infty})$. If $M < \Gdag$, then $U_{\infty}$ is not contained in 
$(\Gdag)^{\prime}$ since $U_{\infty}/(U_{\infty}\cap M)$ is nilpotent. 
But since $H$ acts without fixed point on $U_{\infty}$, we have $U_{\infty} =[U_{\infty},H] \leq 
(\Gdag)^{\prime}$. It follows $\Gdag =M$. 
\end{bew}

In the following two cases we are immediately able to determine the Moufang set.

\begin{prop}\label{H even}
 If $H$ has even order, then $M(U,\tau) \cong 
M({\mathbb F}_q)$.
\end{prop}
\begin{bew} If $H$ has even order, then $H$ contains 
an involution $t$. Since $t$ has no fixed points on $U$, $t$ must 
invert every element in $U$. This implies that $U$ is abelian and hence 
$M(U,\tau) \cong M({\mathbb F}_q).$ 
\end{bew}

\begin{prop} If $|U| \equiv \ 1 \ mod \ 4$, then
$M(U,\tau) \cong M({\mathbb F}_q)$ with $q=|U|$.
\end{prop} 

\begin{bew} If $|U| \equiv \ 1 \ mod \ 4$ and $H$ has odd order, 
then $|\Gdag| = (|U|+1) |U| |H|
\equiv 2 \ mod \ 4$. Hence $\Gdag$ possesses a normal subgroup $L$ of index 
$2$, which contradicts \ref{properties U} (c). 
\end{bew}

From now on, let $N$ be the stabilizer of the set $\{\infty,0\}$ in $\Gdag$. Since 
$H$ has no fixed points apart from $\infty$ and $0$, we have $N=N_{\Gdag}(H)$. For all  
$a \in U^{\#}$ we have $\mu_a \in N \setminus H$, therefore $N =\langle H,\mu_a \rangle$ and 
$|N:H|=2$.

\begin{lemma}\label{H odd}
 If $H$ has odd order, then there is a unique conjugacy class of 
involutions in $\Gdag$.
\end{lemma} 

\begin{bew} Let $s$ and $t$ be involutions in $\Gdag$. Since $\Gdag$ acts $2$-transitively 
on $X$, there are involutions in $N$ which are conjugate to  $s$ and $t$, respectively. Therefore, we may assume that $s$ and $t$ are in $N$.
 Since $|N:H|=2$ and neither $s$ nor $t$ is in $H$, 
the element $st$ must be in $H$ and therefore has odd order. Thus $s$ and $t$ are conjugate.
\end{bew}

\subsection{Zassenhaus Moufang sets with $|U|$ even}

From now on, we assume that the order of $U$ is even. Notice that since $H$ acts regularly on $U$, the order of $H$ is odd.

\begin{lemma}\label{description U}
\begin{enumerate}
\item There is a single $H$-orbit of involutions in $U$. This orbit  is contained in $Z(U)$.
\item For every $a\in U^{\#}$ there is exactly one involution $b \in U$ 
with $\mu_a=\mu_b$. Especially, $\mu_a$ is always an involution.
\item $H$ is cyclic and every element in $N\setminus{H}$ inverts every element in $H$.
\end{enumerate}
\end{lemma} 

\begin{bew} 
(a) All involutions in $U$ are $H$- conjugate by \ref{order mua mub}. Since $U$ is a finite nilpotent group, they are 
all contained in $Z(U)$.

(b) For $x \in U$ the element $\mu_x$ is contained in $N \setminus H$. Suppose $x \in U$ is an involution. 
Then $\mu_x$ is an involution as well. Since $H$ acts regularly on the set of involutions in $U$ 
and since $\mu_x^h = \mu_{xh}$ for all $h \in H$, 
we have by \ref{mua = mub} that $|\mu_x ^H| =|H| = |N \setminus H|$. Thus the claim follows.

(c) If $a$ is an involution, then \ref{mua = mub} implies 
that $C_H(\mu_a)=1$. Since $\mu_a$ is an involution, $\mu_a$ inverts every element in  $H$.
Hence $H$ is abelian.
Since $H$ acts freely on $U$, it follows that $H$ is  cyclic. 
\end{bew}

Since all involutions in $U$ are $H$-conjugate, either all involutions 
are special or no involution is special. We first treat the case 
that all involutions are special and show that this implies $M(U,\tau)$ special 
and hence $\Gdag \cong (P)SL_2(q)$ with $q =|U|$ a power of $2$.
\bigskip\\
\noindent
{\bf Proof of Theorem \ref{theorem1}.}
With \ref{properties U} and \ref{description U} we see that $M(U,\tau)$ satisfies the conditions in \ref{condspecial}. 
Thus $M(U,\tau)$ is special and $U$ is of exponent $2$. By \ref{description U} (c) the Hua subgroup $H$ is abelian, thus
the theorem follows with \cite{G} (alternatively, we can apply \cite{DS2}). 
 \hfill $\qed$

\section{Suzuki Moufang sets}

\subsection{Suzuki $2$-Groups and Suzuki Moufang sets}

\begin{definition}\label{Suz 2 group} [\cite{HB2} VIII, 7.1] A finite group $G$ is called 
a Suzuki $2$-group if the following hold:
\begin{enumerate}
\item $G$ is a nonabelian $2$-group.
\item $G$ has more than one involution.
\item There is a soluble subgroup of $\Aut(G)$ which permutes the involutions 
transitively.
\end{enumerate}
\end{definition}

In the begining of VIII, 7 in \cite{HB2} it is shown that  in fact we may  assume 
\begin{description}
\item[(c')] there is a cyclic subgroup of $\Aut(G)$ which permutes the involutions  transitively. 
\end{description}

We already know 
by \ref{description U} that a Sylow $2$-group of the root group $U$ of a Zassenhaus Moufang set of finite even order is 
either abelian or satisfies (a), (b) and (c').

\begin{example}\label{Beispiel Suz 2 group}[see \cite{HB2}, VIII, 6.7, and \cite{VM}, 2.2]
Let $K$ be a field of characteristic $2$ and let 
$\theta $ be a non-zero endomorphism of $K$. Let $k$ denote the image of $K$ under $\theta$. Let $L$ be a $k$-subvectorspace 
of $K$ with $1 \in L$. (We may assume that $K=k[L]$). Set 
$$A(K,L,\theta):=L \times L$$
 with addition
 $$(a,b) + (c,d) := (a+c, b+d + ac^{\theta}).$$
Then $(A(K,L,\theta),+)$ is a group with 
neutral element $(0,0)$ and inverse element $$-(a,b) =(a,b+a^{1+\theta}).$$

If $K={\mathbb F}_{2^n}$ is finite then $L=K=k$ and $\theta$ is an automorphism of $K$.  We denote this 
group by $A(n,\theta)$.

 If $\theta$ is not the identity, then $A(K,L,\theta)$ is non-abelian and 
the center consists of all elements with first coordinate zero. 
These are exactly the elements of order at most $2$. 

If $\theta$ is a 
{\em Tits endomorphism}, this means $\theta^2$ is the Frobenius endomorphism, 
then $[(a^{-1},0),(a^{\theta},0)] =(0,a+1)$ which implies that  the center equals the 
derived subgroup.

Suppose that $\theta \ne 1$ and that  
$\lambda \mapsto \lambda^{1+\theta}$ is a bijection of $L$ (this is true if $\theta$ is a Tits endomorphism since 
$$\lambda^{(1+\theta) (\theta-1)} = {{\lambda^{2+\theta}} \over {\lambda^{1+\theta}}}   = \lambda;$$
 in the finite case 
$1 + \theta$ is bijective iff the order of $\theta$ is odd). Since $(a,b) \cdot 2 = (0,a^{1+\theta})$, it follows
that the map $x \mapsto x \cdot 2$ induces a bijection from $A(K,L,\theta)/Z(A(K,L,\theta))$ onto
$Z(A(K,L,\theta))$.

Moreover, 
for every $\lambda \in K^{\#}_L:=\{\lambda \in K^{\#}; \lambda L = L\}$ the map 
$$h_{\lambda}: A(K,L,\theta) \to A(K,L,\theta): (a,b) \mapsto
(\lambda a, \lambda^{1+\theta} b)$$
is an automorphism of $A(K,L,\theta)$ and the map $\lambda \mapsto h_{\lambda}$ is an 
injective homomorphism from $K^{\#}_L$ 
into $\Aut(A(K,L,\theta))$ whose image we denote $\Lambda$. 

Suppose that $K={\mathbb F}_{2^n}$ and that the order of $\theta$ is odd.
Then the cyclic group $\Lambda$ acts regularly on the set of involutions of $A(K,L,\theta)$. 
This shows that  $A(n,\theta)$ is a Suzuki $2$-group, see Definition~\ref{Suz 2 group}. 
\end{example}

The finite Suzuki $2$-groups have been classified by Higman \cite{Hi}. His proof uses basically only
linear algebra (see also \cite{HB2}, Theorem VIII, 7.9).

\begin{theorem}\label{suz 2 group}\cite{Hi}
 If $G$ is a finite Suzuki $2$-group, then the exponent of $G$ is $4$, all elements in $Z(G)$ have order at most $2$ 
and either $G \cong A(n,\theta)$ with $o(\theta)$ odd or $|G| =|Z(G)|^3$.
\end{theorem}

Moreover Suzuki showed the following basic fact:

\begin{lemma}\label{Aut suz 2 group}[\cite{Suz1}, Lemma 6]
 If $U=A(n,\theta)$ with $o(\theta)$ odd and if
$H \leq Aut( U)$ is cyclic of order $2^n-1$ and acts transitively on $Z(U)^\#$, then
there is an automorphism $\varphi$ of $U$ such that $\varphi^{-1} H \varphi = \Lambda$.
\end{lemma}

Notice if $H$ is a cyclic subgroup of $Aut( U)$ of order $2^n-1$, then it already acts transitively on $Z(U)^\#$,
as $|U| = 2^n$ and as  an element
of $H$ which centralizes $Z(U)$, already centralizes $U$. 

\begin{definition}\label{Def Suz Moufang}[see \cite{VM}, 2.2] ((Generalized) Suzuki Moufang sets)
Let $K, \theta, k$ and $L$ be as above. Moreover, suppose that $\theta$ is a Tits 
endomorphism. If $K$ is a finite field of order $2^n$, then this implies that $n$  is odd.
Note that $K^2 \subseteq k$ and that $a^{-1} = a^{-2} a \in L$ for all
$a \in L^{\#}$.  

For $a,b \in K$ set $$N(a,b) :=a^{2+\theta} +ab + b^{\theta} .$$

Since $$N(a,b) =\left( {b \over a}\right)^{1 +\theta} + \left(a^{\theta} +{b \over a}\right)^ {1+\theta}
~\mbox{for}~ a \ne 0,$$
it follows that $N(a,b) = 0$ implies $a=b=0$. Set $$U:=A(K,L,\theta).$$
Let $\tau$ be 
the permutation on $U^{\#}$ defined by 
$$(a,b)\tau = 
\left({ b \over {N(a,b)}}, {a \over {N(a,b)}}\right).$$ 

In [20, 2.2] it is shown that $M(U,\tau)$ defines a Moufang set. We call it 
$$\MSuz(K,L,\theta)$$ 
or for $K={\mathbb F}_{2^n}$ just $\MSuz(2^n)$. These Moufang sets are also called (generalized) Suzuki
Moufang sets. The little projective group corresponding to such a Moufang set is called generalized Suzuki 
group or $\Suz(K,L,\theta)$ and $\Suz(2^n)$ if $K={\mathbb F}_{2^n}$.

An easy but tedious computation shows that $\tau^2 =1$ and that $\tau \mu_{(a,b)}$ induces the automorphism 
$$h_{N(a,b)^{2-\theta}}$$ on $U$ (see 2.2 in \cite{VM} where a matrix representation for $\Suz(K,L,\theta)$ is given).
Especially since $H \leq \Lambda$ and $\Lambda$ acts freely on $U$, it follows that $\MSuz(K,L,\theta)$ is a Zassenhaus Moufang set and $$\mu_{(a,b)} = \mu_{(c,d)}~\mbox{iff}~N(a,b) =N(c,d).$$
\end{definition}

\subsection{The Case: no special involution}

From now on we assume that no involution of $U$ is special. We are going to prove that $M(U,\tau)$ is a 
Suzuki Moufang set. 

\subsubsection{The Identification of $U$}

First, we need that $U$ is a $2$-group.

\begin{theorem}\label{p-group}
If $M(U,\tau)$ is a Zassenhaus Moufang set of arbitrary finite order, then $U$ is a $p$-group
 \end{theorem}
\begin{bew}
This was proven by Feit \cite{F}. His proof with some improvement by Bender is contained in 
\cite{HB3}. More precisely the assertion follows from 4.1, 6.3, 6.5, 6.6 and 5.7 of \cite{HB3}.
\end{bew}

By \ref{description U} $U$ is either an abelian $2$-group or a Suzuki $2$-group.
 If $U$ is abelian, then $M(U,\tau)$ is 
special by \cite{S}. As we assume that $U$ does not contain special involutions, it follows 
that  $U$ is a Suzuki $2$-group. By 3.9 (b) and by 5.4 this means that $U$ does not contain special elements at all.

From now on, let $q:=|Z(U)|$. Then $q$ is a power of $2$ and $|H|=q-1$. 
Moreover, $|U| =q^2$ or $|U|=q^3$, and $Z(U)$ has exponent $2$ by \ref{suz 2 group}. Since $M(U,\tau)$ is proper, $q>2$.

If $a \in Z(U)^\#$, hence $a$ is an involution, then $\mu_a$ is an involution 
and  $$\alpha_a^{\mu_a} =\mu_a^{\alpha_{-\sim a}}$$ by \ref {a inv} (a). The last equals Suzuki's
structure equation (XI, 10.6 in \cite{HB3}).

Since $a$ is not special, we have $\sim a \neq a$, see \ref{equivalence special}.
We first study the subgroup $U$ before we show that $M(U,\tau)$ is a Suzuki Moufang set.

\begin{lemma}\label{sim a 2} For all involutions $a \in U$, we have $(\sim a) \cdot 2 =a $.
\end{lemma}

 \begin{bew} Set $D := \langle \alpha_a, \mu_a \rangle$. Then $D$ is
 dihedral 
since $\alpha_a$ and $\mu_a$ are involutions. The order of $\alpha_a \mu_a$ is 
odd since $$C_{\Gdag}(\alpha_a) \cap C_{\Gdag}(\mu_a) = U_{\infty} \cap
 C_{\Gdag}(\mu_a) =1.$$
Set $E :=\langle \mu_a \alpha_a \rangle$. Then $\alpha_{\sim a} \in N_G(E)$,
since by \ref{a inv} (a) we have that $$D^{\alpha_{-\sim a}} = \langle \alpha_a^{\alpha_{-\sim a}}, 
\mu_a^{\alpha_{-\sim a}}\rangle = \langle \alpha_a, \alpha_a^{\mu_a} \rangle=\langle \alpha_a, \mu_a \rangle =D$$
and since $E = D^{\prime}$ is characteristic in $D$. On the other hand, 
$C_{U_{\infty}}(E) =C_{U_{\infty}}(\mu_a) = 1$. Hence $N_{U_{\infty}}(E) \leq \Aut(E)$
is abelian. If $\alpha_t \in N_{U_{\infty}}(E)$ is an involution with $t\ne a$, then $\alpha_t$ or 
$\alpha_{a+t}$ fixes a point in $E$. Since $C_{\Gdag}(\alpha_x) \leq U_{\infty}$ for all $x \in U^{\#}$ this 
implies $U_{\infty} \cap E \ne 1$. But $U_{\infty} \cap E$ is centralized by $\alpha_a$ and thus also by 
$\mu_a$, a contradiction since $U_{\infty} \cap C_{\Gdag}(\mu_a) =1$. 
It follows that $\alpha_a$ is the unique involution in 
$N_{U_{\infty}}(E)$. So the latter group is cyclic. Since $U_{\infty}$ has exponent $4$, 
the claim follows. 
\end{bew}

Next we show that $|U|=q^2$ and that therefore $U\cong A(n,\theta)$ for some $n \in \NN$.

 \begin{lemma}\label{order 3}
$H$ does not contain an element of order $3$. 
\end{lemma}

\begin{bew}
Suppose otherwise. Then there are involutions $a,b \in U$ 
such that $\mu_a \mu_b$ has order $3$. Thus also $$(\mu_a \mu_b)^{\alpha_{-\sim a} \mu_a} =
\mu_a^{\alpha_{-\sim a} \mu_a} \mu_b^{\alpha_{-\sim a} \mu_b} = \alpha_a \mu_b^{\alpha_{-\sim a} \mu_a} $$
has order $3$ (we obtained the second equality by applying \ref{a inv} (a)). Since $\sim b$ is the unique fixed point of $\mu_b$, we have that $$x:=(\sim b -\sim a) \mu_a = (\sim b)^{\alpha_{-\sim a} \mu_a}$$ is the 
unique fixed point of $\mu_b^{\alpha_{-\sim a} \mu_a}$. Since $a \ne b$, we have $x \ne \infty$. 
Thus for $t:=(\mu_b^{\alpha_{-\sim a} \mu_a})^ {\alpha_{-x}}$
we see, as $a \in Z(U)$, that $\alpha_a^{\alpha_{-x}} t = \alpha_a t$ has order $3$ and $t$ fixes $0$. 
Since all involutions of $\Gdag$ which fix $0$ lie in $U_0$, there is an involution $c \in U$ 
such that $t = \alpha_c^{\mu_a}$. 

Thus $1 = (\alpha_a t)^3 = (\alpha_a \alpha_c^{\mu_a})^3 $, which implies that $1 = (\alpha_a^{\mu_a} \alpha_c)^3 $.
Hence with \ref{a inv} (a) we obtain $1 = (\mu_a^{\alpha_{-\sim a}}\alpha_c)^3$. As $c$ is in $Z(U)$, we have
$\mu_a^{\alpha_{-\sim a}}\alpha_c = (\mu_a \alpha_c)^{\alpha_{-\sim a}}$. This implies that $(\mu_a \alpha_c)^3 = 1$.

By multiplying the last equation with $\mu_a$ from the right we get
$$\mu_a = \alpha_c^{\mu_a} \alpha_c \alpha_c^{\mu_a}.$$
Since $\mu_c$ is the unique element in $U_0 \alpha_c U_0$ interchanging 
$0$ and $\infty$, this implies $\mu_a =\mu_c$ and thus $a=c$. Now by \ref{cond special} this equation implies that $a$ 
is special, a contradiction to our assumption that there is no special element. 
\end{bew}

In fact we showed that if $H$ contains an element of order $3$, then $Z(U)$  contains a special element
and therefore $M(U,\tau) = M(q)$ by Theorem~\ref{theorem1}.

\begin{lemma}\label{Antheta} 
\begin{enumerate}
\item $q=2^n$ with $n$ odd. 
\item Let $a \in U$ be an involution. Then the order of $\mu_a \alpha_a$ is $5$.
\item $U$ is isomorphic to $A(n,\theta)$.
\end{enumerate}
\end{lemma}
\begin{bew} We have taken this proof from \cite{HB3}, XI, 11.2.

(a) If $n$ was even, $3$ would divide $q-1$. This is not possible because of \ref{order 3}.

(b) We compute using \ref{a inv} (a) in the second equality below that $$(\alpha_a \mu_a)^{\alpha_{-\sim a}} = \alpha_a \mu_a^{\alpha_{-\sim a}} =\alpha_a \alpha_a^{\mu_a}=(\alpha_a \mu_a)^2.$$ The next equality follows 
as $\alpha_a$ is an involution and the second because of the last calculation and as ${(\sim a) \cdot 2} =a$ by \ref{sim a 2}:
 $(\alpha_a \mu_a)^{-1} =
(\alpha_a \mu_a)^{\alpha_a} = (\alpha_a \mu_a)^4$ and thus $(\alpha_a \mu_a)^5 =1$.

(c) If $U$ is not isomorphic to $A(n,\theta)$, then $|U|=q^3$ and thus $$|\Gdag| ={(q^3+1)q^3(q-1)}~\mbox{ with}~q=2^n,~
n~\mbox{ odd}.$$
Since $5$ divides $|\Gdag|$ but neither $q^3$ nor $q-1$, $5$ divides $q^3 +1$ and thus $5$ divides ${q^6-1} ={2^{6n}-1}$. This implies
$2^{6n} \equiv \ 1 \ mod \ 5$ and hence $4|6n$, a contradiction to $n$ odd.
 \end{bew}

Thus $|U|=q^2$ and $|\Gdag| =(q^2+1)q^2(q-1)$, which is exactly the order of the Suzuki group $\Suz(q)$.

To identify $G$ as $\Suz(q)$ we still need to find an involution $\tau$ in $G$ which permutes the elements of $U^\#$ as described in Definition~\ref{Def Suz Moufang}. Our next aim is to find a convenient partition of $U$, see \ref{partition}, which will help us to calculate the action of a possible $\tau$ on $U^\#$.

\subsubsection{A partition  of $U$}

\begin{lemma}\label{represent} 
The set $\{0\} \cup \{\sim a~|~ a \in Z(U)^{\#}\}$ is a system 
of representatives for $U/Z(U)$. 
\end{lemma}

\begin{bew} By \ref{Antheta}, $U \cong A(n,\theta)$. Thus if $|Z(U)|=q$, 
then $q =|U:Z(U)|$. By \ref{sim a 2}, $(\sim a) \cdot 2 =a$ for all $a \in Z(U)^{\#}$.
Thus the claim follows since by \ref{Beispiel Suz 2 group} a subset $X$ of $U$ is a system of representatives for $U/Z(U)$ if and only 
if the map $x \mapsto x \cdot 2$ induces a bijection between $X$ and $Z(U)$.
\end{bew}

\begin{lemma}\label{orbits2}
If $a,b$ are two different involutions in $Z(U)$, then 
$\mu_{a+b} =\mu_{-\sim a+\sim b}$.
\end{lemma}

\begin{bew} This follows from \ref{orbit of subgroup of H} with $V=Z(U)$ and $K=H$.
\end{bew}

The sets $Z(U)^{\#}$, $\{\sim a~|~a \in Z(U)^{\#}\}$ and $\{-\sim a~|~a \in Z(U)^{\#}\}$ are three
orbits of $H$ on $U$ (so each is of size $q-1$).
We will show that each of the remaining $(q^2 -1) - 3(q-1) =(q-2)(q-1)$ non-trivial elements of $U$ can uniquely be written as the sum $-\sim a +\sim b$ with $a,b 
\in Z(U)^{\#}$ and $ a \ne b$. Having succeeded we will know the action of $\mu_a$ on $U^\#$ for all $a \in U$.

The next very technical lemma is needed for the proof of the very important Lemma~\ref{-mua + mub}.

\begin{lemma}\label{sim} 
\begin{enumerate}
\item If $a,b \in Z(U)^{\#}$ are different, then $$\sim (-\sim a +\sim b) = -\sim (a\mu_b \mu_{a\mu_a \mu_b +b}) +\sim({a\mu_a \mu_b +b})$$
$$= -\sim (a +a\mu_a \mu_b) +\sim (b +a\mu_a \mu_b).$$
\item
If $e \in Z(U)^{\#}$ and $g,h \in H$ are different, then
$$\sim(-\sim eg + \sim eh) =-\sim (eh^2 g^{-1} + eg) +\sim (eh^2 g^{-1}+eh). $$

\end{enumerate} 
\end{lemma}

\begin{bew} 
(a) By \ref{mu and sim} (e)
$$\sim (-\sim a +\sim b) = (-(-\sim a +\sim b)\mu_b)\mu_b =
(-(a\mu_a \mu_b^2 -b\mu_b)\mu_b)\mu_b = $$
$$(-((a\mu_a \mu_b -b)\mu_b +\sim b))\mu_b = 
(b \mu_b - (a\mu_a \mu_b +b)\mu_b) \mu_b =$$
$$a\mu_a \mu_b \mu_{a \mu_a \mu_b +b} + \sim (a\mu_a \mu_b +b)
=(-\sim a) \mu_b \mu_{a\mu_a \mu_b +b} + \sim (a\mu_a \mu_b +b)$$
$$=-\sim (a\mu_b \mu_{a\mu_a \mu_b +b}) +\sim({a\mu_a \mu_b +b}).$$
Now with \ref{orbits2}
$$\mu_{a+b} = \mu_{-\sim a +\sim b} =\mu_{\sim(-\sim a +\sim b)} =\mu_{-\sim (a\mu_b \mu_{a\mu_a \mu_b +b}) +\sim({a\mu_a \mu_b +b})}$$
$$=\mu_{a\mu_b \mu_{a\mu_a \mu_b +b} +a\mu_a \mu_b +b}$$
and thus
$a +b =a \mu_b \mu_{a\mu_a \mu_b + b} +a\mu_a \mu_b + b$. 
This finally implies $a\mu_b \mu_{a\mu_a \mu_b +b} = a +a\mu_a \mu_b$.

(b)
This follows from (a) with $a = eg$ and $ b = eh$. One has  
$$\mu_a \mu_b = \mu_{eg} \mu_{eh} = g^{-1} \mu_e gh^{-1} \mu_e h =
g^{-2} h \mu_e^2 h = g^{-2} h^2$$ 
and therefore $a\mu_a \mu_b = eg g^{-2} h^2 = eg^{-1} h^2$.
\end{bew}

\begin{lemma}\label{-mua + mub}
 If $a,b,c \in Z(U)$ are different involutions, then
$-\sim a +\sim b \not \in \{c,-\sim c, \sim c\}$. 
\end{lemma}

\begin{bew} By \ref{represent} the elements $-\sim a$ and $-\sim b$ lie in different cosets of 
$U/Z(U)$, hence $-\sim a +\sim b$ cannot be in $Z(U)$. If $-\sim a+\sim b
= -\sim c$ for some $c \in Z(U)$, then $-\sim b +\sim a =\sim c$, so we only have to show 
that $-\sim a +\sim b = \sim c$ cannot hold. 

Suppose that $\sim c =-\sim a +\sim b$ with $a\ne b$.
Then with \ref{sim} (a) we get $c =\sim (-\sim a +\sim b) = -\sim (a +a\mu_a\mu_b) +\sim (b +a\mu_a \mu_b)$. 
But we have already shown above that this is impossible.
\end{bew}

From now on, $e$ will be a fixed involution in $U$ and 
$$\tau:=\mu_e.$$

We introduce coordinates.
By \ref{Antheta}, $U$ is isomorphic to $A(n,\theta)$ with $o(\theta)$ odd.
We can label the elements as $(a,b)$ with $a,b \in K := {\mathbb F}_{2^n}$.

We claim  that we can assume that 
$H= \Lambda$ with $\Lambda \leq \Aut(U)$ as defined in \ref{Beispiel Suz 2 group}.
By \ref{Aut suz 2 group} there is an automorphismus $\varphi$ of $U$ 
such that $\varphi^{-1} H\varphi = \Lambda$.
So if $h \in H$ then there is a $\lambda \in K$ with $h^{\varphi} = h_{\lambda}$. 
Then $$(a,b)\varphi^{-1} h  = (a,b) h_{\lambda} \varphi^{-1} = (a \lambda, b \lambda^{1 +\theta}) \varphi^{-1}.$$
Hence writing $(a,b)$ instead of $(a,b)\varphi^{-1}$ we can assume $H = \Lambda$. 

Since $H$ acts transitively on the set of involutions, we can also assume $e =(0,1)$.

\begin{lemma}\label{commutator} 
\begin{itemize}
 \item[(a)] Let $K_\theta = Fix_K(\theta)$. Then $K_\theta = \FF_2$.
\item[(b)] $[\sim e, \sim eh] \ne 1$ for $h$ in $H^\#$.
\end{itemize}
\end{lemma}

\begin{bew} 
In $A(n,\theta)$, an element of the form $(x,y)$ with $x \in K_\theta^{\#}$ commutes 
with $(x^{\prime}, y^ {\prime})$ if 
and only if $x^{\prime} \in K_\theta$, see Example~\ref{Beispiel Suz 2 group}. Let $a=(0,u)$ and $b=(0,v)$ be involutions. 
Since $(\sim a )\cdot 2 =a$ and since $ (\sim b )\cdot 2 = b$ by \ref{sim a 2}, we get $\sim a = (u^{(1+\theta)^{-1}},x)$ and 
$\sim b = (v^{(1+\theta)^{-1}},y)$ with 
$x,y \in K$.  Clearly $\lambda \mapsto  \lambda^{1+\theta}= \lambda^2$ induces a bijection of $K_\theta$.
Thus if $u \in K_\theta$, then $\sim a$ and $\sim b$ commute iff $v \in K_\theta$.

 Suppose that $|K_\theta| >2$ and that $\lambda \in K_\theta \setminus \{0,1\}$. Set 
$h:=h_{\lambda}$, see Example~\ref{Beispiel Suz 2 group}.
Then $eh = (0,\lambda^{1+\theta}) = (0,\lambda^2)$, hence $\sim e$ and 
$\sim eh$ commute. Therefore
$$(-\sim eh  +\sim e) \cdot 2 =(-\sim eh)\cdot 2 + (\sim e)\cdot 2 = eh +e = 
(\sim (eh +e)) \cdot 2,$$ hence $c:=-\sim eh +\sim e - \sim (eh +e) \in Z(U)$, see Example~\ref{Beispiel Suz 2 group}. 
Note that $c \ne 0$ by \ref{-mua + mub}.
In the next calculation the first equality follows from \ref{mu and sim} (k) and the third from \ref{sim} (b).
 $$(-\sim e +\sim eh)\mu_{e +eh} = -\sim -(-\sim e +\sim eh ) =
-\sim (-\sim eh +\sim e) = $$
$$-(-\sim (eh + eh^{-1}) +\sim (e + eh^{-1}) ) = -\sim (e +eh^{-1}) + \sim (eh +eh^{-1}).$$ 
As $eh + eh^{-1} = (0,\lambda^2) + (0,\lambda^{-2})= (0,\lambda^2 +\lambda^{-2})$ and $e + eh^{-1} = (0,1) + (0,\lambda^{-2}) =
(0,\lambda^4 +1)$, it follows that $\sim (e +eh^{-1}) $ and $\sim (eh + eh^{-1})$ commute.
Then by the last calculation and by \ref{sim a 2} $$(-\sim e +\sim eh)\mu_{e +eh}\cdot 2 = (-\sim (e +eh^{-1}))\cdot 2 + (\sim (eh +eh^{-1}))\cdot 2
=$$ $$ (e +eh^{-1}) + (eh +eh^{-1}) = eh + e.$$
Therefore as $eh + e \in Z(U)$ and as $(- \sim (e + eh))\cdot 2 = eh + e$ as well, we get  
$${d:=-\sim (eh^{-1} + e)} + {\sim (eh + eh^{-1}) - \sim (e + eh)} \in Z(U).$$
 Again by \ref{-mua + mub}, 
$d\ne 0$.
Now by  our calculation above $${d = (-\sim eh + \sim e)\mu_{e+eh}} - \sim (e + eh)$$ and 
$$c=(-\sim  eh +\sim e) - {\sim (e +eh)}.$$
With \ref{orbits2} and \ref{mu and sim} (l) applied to $a=(-\sim eh + \sim e)$ and $b = \sim (e + eh)$ we get $\mu_c = \mu_d$. 
But since $c$ and $d$ are involutions, \ref{mua = mub} implies $c=d$. Therefore,
$$(-\sim eh + \sim e)\mu_{e+eh}  = (-\sim  eh +\sim e)) .$$
With \ref{a inv} (b) we get $-\sim  eh +\sim e = \sim (e + eh)$. But this contradicts \ref{-mua + mub}. Hence $K_{\theta}={\mathbb F}_2$ 
and therefore by the second paragraph of the proof
$C_U(\sim e) = \langle Z(U), \sim e\rangle$. 
\end{bew}

The next two lemmata explicitly describe the elements $\sim a, -\sim a$ and
$\sim a + \sim b$, where $a \neq b$ are involutions in $U$.

\begin{lemma} \label{formula}
After possibly replacing $\theta$ by $\theta^{-1}$, we can 
assume that $$-\sim e = (1,0).$$
Then we have for all $a \in K^{\#}$
$$-\sim (0,a^{1+\theta}) =
(a,0)~\mbox{and}~\sim (0,a^{1+\theta}) =(a,a^{1+\theta}).$$
\end{lemma}
\begin{bew}
By \ref{commutator} the automorphism $\theta$ has no fixed points other than $0$ and $1$.
Thus we can prove the first statement as in Lemma XI, 11.12 in \cite{HB3}. 
The second holds since $(0,a^{1+\theta}) = e h_a$ and hence 
$-\sim (0,a^{1+\theta}) = -\sim eh_a = (-\sim e)h_a = (a,0)$ by \ref{mu and sim} (h).
\end{bew}

\begin{lemma}\label{formula1}
Let $(a,b) \in U$ with $0 \not\in \{a,b\}$ and $b \ne 
a^{1+\theta}$. Set $t = \left({b \over a}\right)^{\theta^{-1}}$ and 
$s =a-t$.
\begin{enumerate}
\item 
$(a,b) =-\sim(0,s^{1+\theta}) +\sim (0,t^{1+\theta})$.
\item
If $(a,b) = -\sim (0,u^{1+\theta} ) + \sim (0,v^{1+\theta})$, then $s=u$ and $t=v$.
\end{enumerate}
\end{lemma}

\begin{bew} 
By \ref{formula} the formula in (a) holds iff $(a,b) = (s,0) -(t,0) = (s,0) +(t,t^{1+\theta})
=(s+t, st^{\theta} + t^{1+\theta})$. 
Now $s+t =a$ and $t^{1+\theta} + st^{\theta} = t^{\theta} (s+t) = a \cdot {b \over a} = b$,
which implies (a).

By (a) all elements in $U$ which are not in $Z(U)~ \cup \sim Z(U)^{\#} \cup -\sim Z(U)^{\#}$ 
are of the form $-\sim x +\sim y$ with $x,y \in Z(U)^{\#}, x \ne y$. Since we have 
$${|U \setminus (Z(U) \cup \sim Z(U)^{\#} \cup -\sim Z(U)^
{\#}) |} $$ $$= q^2 -q - 2(q-1) = (q-1)(q-2),$$ we see that $x$ and $y$ are always uniquely determined. 
This shows (b).
\end{bew} 

As a consequence of  \ref{formula1} we obtain:

\begin{coro} \label{partition}
The following is a partition of $U$:
 $$U = Z(U)~ \cup ~\sim Z(U)^{\#}~ \cup ~ -\sim Z(U)^{\#}~
\cup ~\{-\sim a+\sim b~|~ a,b \in Z(U)^{\#}, a\ne b\}.$$
\end{coro}

We will call such a partition, which depends of the Moufang set $M(U,\tau)$ and which is
fundamental for the rest of the proof, a {\em Suzuki partition}.

The next lemma is proven to verify the last paragraph of Definition~\ref{Def Suz Moufang}.

\begin{lemma} For $(a,b) \in U^{\#}$, set $N_0(a,b) := 
a^{1+\theta} + a^{\theta-\theta^{-1} } b^{\theta^{-1}} + b$. 
Then $\mu_{(0, N_0 (a,b))} =\mu_{(a,b)}$.
 \end{lemma}

\begin{bew} 
For $a=0$, we have $N_0(a,b) =b$. If $a \ne 0, b =0$ or
 $b=a^{1+\theta}$, then 
$N_0(a,b) =a^{1+\theta}$ and the claim is true by \ref{orbits2}, since by \ref{formula} $-\sim (0,a^{1+\theta}) =
 (a,0)$ and 
$\sim (0,a^{1+\theta}) = (a,a^{1+\theta})$. If $a \ne 0$ and $ b$ neither $0$ nor $a^{1+\theta}$,
 then 
$N_0(a,b) = s^{1+\theta} + t^{1+\theta}$ with $s,t$ as is \ref{formula1}, and so this case
 follows again from  \ref{orbits2}.
\end{bew}

\subsubsection{The action of $\tau$ on $U$}

In the next two lemmata we calculate the action of $\tau = \mu_e$ on $U$.

\begin{lemma}\label{formula2} 
Let $g,h \in H$ with $g\ne h$. Then 
$$(-\sim eg +\sim eh)\tau =
-\sim ej^{-1}h^{-2} +\sim eh^{-1},$$ where $j \in H$ such that 
$ej = eg^{-1} + eh^{-1}$. 
\end{lemma}

\begin{bew} Let $j \in H$ with $ej = eg^{-1} + eh^{-1}$. Then, as $\tau$ inverts every element in $H$, \ref{description U} (c), we get with (3A)
and \ref{mu and sim} (g) that 
$$(-\sim eg + \sim eh) \tau =
(e\tau g -e\tau h)\tau = (eg^{-1} \tau - eh^{-1} \tau) \tau = $$
$$(eg^{-1} +eh^{-1})\mu_{eh^{-1}} + \sim eh^{-1} = 
(eg^{-1} +eh^{-1}) h \tau h^{-1} + \sim eh^{-1} = $$
$$ej \tau h^{-2}  + \sim eh^{-1}
= e\tau j^{-1}h^{-2} +\sim eh^{-1}= -\sim ej^{-1}h^{-2} + \sim eh^{-1}.$$ 
\end{bew}

\begin{lemma}\label{tau}
If $a,b \in K^{\#}$ with $b\ne a^{1+\theta}$, 
then
\begin{enumerate} 
\item $(0,a^{1+\theta})\tau = (a^{-1},0)$.
\item $(a,0)\tau =(0,a^{-1-\theta})$.
\item $(a,a^{1+\theta})\tau = (a^{-1}, a^{-1-\theta})$.
\item $(a,b)\tau = ({s \over {Nt}} + {1 \over t}, {1 \over {t^{\theta}}}({s 
\over {Nt} } + {1 \over t}))$ with $s,t $ as in \ref{formula1} and $N^{1+\theta} = N_0(a,b)$.
\end{enumerate}
\end{lemma}
\begin{bew} 
As $\tau = \mu_e$, we get with \ref{mu and sim} that
$eh_a\tau = e\tau h_a^{-1} =-\sim e h_a^{-1}$, which shows  (a).
Similarly, $(-\sim e )h_a \tau = (-\sim e )\tau h_a^{-1} = eh_a^{-1}$, which is (b).

By \ref{a inv}(b) 
 $(\sim e) h_a \tau = (\sim e) \tau h_a^{-1} = (\sim e) h_a^{-1}$, which is (c).

Let $s,t$ be as in \ref{formula1}. Set $g:=h_s$ and $h:=h_t$.
Then $(a,b) = -\sim eg + \sim eh$ and hence by \ref{formula2} 
$(a,b)\tau = -\sim ej^{-1} h^2 + \sim eh^{-1}$ with $ej = eg^{-1} +eh^{-1} =(0,s^{-1-\theta}) +(0,t^{-1-\theta}) =
(0,(st)^{-1-\theta}(s^{1+\theta} +t^{1+\theta}) =(0,(st)^{-1-\theta} N_0(a,b))$. 
Therefore $ej^{-1} = (0,(st)^{1+\theta} N_0(a,b)^{-1})$ and so $ej^{-1} h^{-2} = 
(0,s^{1+\theta} t^{-1-\theta} N_0(a,b)^{-1})$ and $-\sim (ej^{-1} h^{-2})= 
(st^{-1} N^{-1},0)$. We get
$$(a,b)\tau = (st^{-1} N^ {-1},0) + (t^{-1}, t^{-1-\theta}) = 
\left({s \over {Nt} } + {1 \over t}, {1 \over {t^{\theta}}}\left({s \over {Nt} } + {1
 \over t}\right)\right) ,$$ showing (d).
\end{bew}

\begin{lemma}\label{tits}
$\theta$ is a Tits automorphism.
 \end{lemma}
\begin{bew}
The strategy of the proof is to apply $\tau$ to a convient element of $U$ which allows us to calculate
$a^\theta$ for $a \in K$.

For all $h \in H^\#$, we have by \ref{tau}, \ref{a inv} (b) and (3A) $$(\sim eh - \sim e)\tau = ((\sim e h^{-1})
 \tau -(\sim e) \tau)\tau =
(\sim eh^{-1} -\sim e)\tau +e.$$
For $a \in K \setminus \{0,1\}$ and $h=h_a$, this means
$$(*) \hspace{2cm}(a+1,a (a^{\theta} +1))\tau = (a^{-1}+1, a^{-1} (a^{-\theta} +1))\tau +(0,1).$$
By \ref{formula1} there are uniquely determined elements $s,t,u,v \in K^{\#}$ such that  
$${(a+1, a (a^{\theta} +1)) = (s,0) -(t,0)}~\mbox{and}~
(a^{-1}+1,a^{-1} (a^{-\theta} +1)) = 
(u,0) -(v,0).$$ We compute 
$$ t = a^{\theta^{-1}} (a+1) (a^{\theta^{-1}} +1)^{-1},$$
$$ s = a+1 + t = (a+1) (a^{\theta^{-1}} +1)^{-1} (a^{\theta^{-1}}+ (a^{\theta^{-1}} +1)) 
= a^{-\theta^{-1}} t,$$
$$ v =(a^{-\theta^{-1}}+1)^{-1}a^{-\theta^{-1}}(a^{-1}+1) =(a^{\theta^{-1}}+1)^{-1}a^{-1}(a+1)=
 a^{-1-\theta^{-1}} t=a^{-1}s $$ and 
$$ u = (a^{-1} +1) + v = a^{-1} (a+1) +a^{-1} s = a^{-1}( a+1+s) = a^{-1} t.$$

Set $$N := (s^{1+\theta} +t^{1+\theta})^{(1+\theta)^{-1}}~\mbox{and }~
 M := (u^{1+\theta} +v^{1+\theta})^{(1+\theta)^{-1}}.$$
Then $M = a^{-1} N$ and $N^{1+\theta} = (1+a^{-1-{\theta}^{-1}}) t^{1+\theta}$.
By substituting first the calculated expressions in equation $(*)$ and then by applying \ref{tau} we get
$$(sN^{-1} t^{-1} + t^{-1}, sN^{-1} t^{-1-\theta} +t^{-1-\theta}) = (uM^{-1}
 v^{-1} +v^{-1}, uM^{-1} v^{-1-\theta} +v^{-1-\theta} +1).$$
Hence
$$ a^{-\theta^{-1}} N^{-1} + t^{-1} =  a^{1+\theta^{-1}} N^{-1} + a^{1+\theta^{-1}}
 t^{-1}.$$
This implies 
$$(a^{-\theta^{-1}} +a^{1+\theta^{-1}}) t = (a^{1+\theta^{-1}}+1) N.$$
Thus
$$ N = a^{-\theta^{-1} }( a^{2\theta^{-1}+1}+1) (a^{1+\theta^{-1}}+1)^{-1} t.$$
Therefore we get 
$$ (1 +a^{-1-\theta^{-1}}) t^{1+\theta} = N^{1+\theta} = a^{-1-\theta^{-1}}
 (a^{2\theta^{-1} +1} +1 )(a^{2+\theta} +1)
(a^{1+\theta^{-1}}+1)^{-1}(a^{\theta+1}+1)^{-1} t^{1+\theta}$$
and so
$$a^{1+\theta^{-1}} (1+a^{-1-\theta^{-1}})(1+a^{1+\theta^{-1}})(1+a^{1+\theta})
 =
(a^{2 \theta^{-1} +1} +1)(a^{2+\theta} +1).$$
Hence 
$$(1+a^{1+\theta^{-1}})^2 (1+a^{1+\theta}) = a^{3 +2\theta^{-1} +\theta} +1 +
a^{2\theta^{-1} +1} + a^{2+\theta}$$
and
$$ 1 +a^{2+2\theta^{-1} } + a^{1+\theta } + a^{3+2\theta^{-1}+\theta} =
 a^{3 +2\theta^{-1} +\theta} +1 +a^{2\theta^{-1} +1} + a^{2+\theta}.$$
We get
$$ a^{2+2\theta^{-1}} + a^{1+\theta} = a^{2+\theta} + a^{2\theta^{-1} +1}$$
and so 
$$ a(a+1) a^{2\theta^{-1}} = a(a+1)a^{\theta}.$$
Since $a \ne 1$, this implies $a^{2\theta^{-1}} = a^{\theta}$ and hence 
$a^2 = a^{\theta^2}$. 
\end{bew}

Now we are able to prove the second main theorem.
\\\\
\noindent
{\bf Proof of Theorem \ref{theorem2}:} By \ref{Antheta} $U$ is isomorpic to 
$A(n,\theta)$ with $|U| =2^{2n}$, and by \ref{tits} $\theta$ is a Tits automorphism. Therefore it remains to show
 that 
$(a,b)\tau = \left({b \over {N(a,b)}}, {a \over {N(a,b)}}\right)$ with $N(a,b) =
a^{2+\theta} + ab + b^{\theta}=N_0(a,b)^{\theta}$. If $a=0$, then $N(0,b) =
 b^{\theta}$. By 
\ref{tau} $$(0,b)\tau =(b^{-(\theta+1)^{-1}},0) = (b^{-\theta+1},0) =
({b \over {b^{\theta}}},0) = \left({b \over {N(0,b)}}, {0 \over {N(0,b)}}\right).$$

If $b=0$, then $N(a,0) = a^{2 +\theta}$ and by \ref{tau}
$$(a,0) \tau = (0,a^{-1 -\theta}) = (0,{a \over {a^{2+\theta}}}) =
\left({0 \over {N(a,0)}}, {a \over {N(a,0)}}\right).$$

If $b =a^{1+\theta}$, then 
$$N(a,b) = a^{2+\theta} + aa^{1+\theta} +a^{\theta^2+1}=a^{2 +\theta}$$ and
$$(a,b)\tau =(a^{-1},a^{-1 -\theta}) =\left({{a^{1+\theta}} \over {a^{2+\theta}}},
{a \over {a^{2+\theta}}}\right) = \left({b \over {N(a,b)}}, {a \over {N(a,b)}}\right).$$

If $a,b \ne 0$ and $b \ne a^{1+\theta}$, then by \ref{tau} (d)
$$(a,b)\tau = \left({s \over {Nt}} + {1 \over t}, {1 \over {t^{\theta}}} 
({s \over {Nt}} + {1 \over t})\right)$$ with $t = ({b \over a}) ^{\theta^{-1}}$,
$s = a + t $ and $N^{1+\theta} =N_0(a,b) = N(a,b)^{\theta^{-1}}$. 
Since $(1 +\theta)^{-1} =1-\theta $, we have 
$$N =( N(a,b)^{\theta^{-1}})^{\theta -1} = N(a,b)^{1- \theta^{-1}}.$$
We compute 
$$st^{-1} = (a + a^{1-\theta^{-1}}b^{-\theta^{-1}} ) a^{\theta^{-1}}
b^{-\theta^{-1}}  = a^{1+\theta^{-1}} b^{-\theta^{-1}} +1$$ and 
$$N(a,b)^{\theta^{-1}} st^{-1} = (a^{1+\theta^{-1}} b^{-\theta^{-1}} +1)
(a^{1+\theta} +a^{\theta^{-1}} b^{\theta^{-1}} +b) =$$
$$ a^{2+\theta +\theta^{-1}}b^{-\theta^{-1}} +a^{1+2\theta^{-1}}  
+ a^{1+\theta^{-1}} b^{1-\theta^{-1}} + a^{1+\theta}  + a^{\theta^{-1}}
b^{\theta^{-1}} + b =$$
$$a^{\theta^{-1}} b^{-\theta^{-1}}(a^{2+\theta}+ab + b^{2\theta^{-1}}) + b
+ a^{1+\theta} + a^{1+2\theta^{-1}} = {1 \over t} N(a,b) + b.$$
Thus $${s \over {Nt}} + {1 \over t} ={1 \over {N(a,b)}} \left({1 \over t} N(a,b) 
+b\right) + {1 \over t} = {b \over {N(a,b)}}$$ and 
$${1 \over {t^{\theta} }} \left({s \over {Nt}} + {1 \over t}\right) = {a \over b} 
{b \over {N(a,b)}} = {a \over  {N(a,b)}}.$$
\hspace*{1cm} \hfill $\qed$

\section{Generalized Suzuki Moufang sets}

If $M(U,\tau)$ is a finite Suzuki Moufang set, then we have seen that the following hold:
\begin{enumerate} 
\item $H$ is transitive on $Z(U)^{\#}$.
\item For every $a \in U^{\#}$ there is a $b \in Z(U)^{\#}$ with $a \in V_b$.
\item $U$ has a Suzuki partition, this means 
$$U =$$ $$Z(U)~ \cup \sim Z(U)^{\#}~ \cup~ -\sim Z(U)^{\#}~ \cup
~{\{-\sim a +\sim b~|~a,b \in Z(U)^{\#}, a \ne b\}}.$$
\end{enumerate}

In the infinite case, one can generalize the concept of Suzuki Moufang 
sets. It turns out that a generalized Suzuki Moufang set $M(K,L,\theta)$ is an 'ordinary' Suzuki 
Moufang set (this means $\theta$ bijective and hence $K=L=K^{\theta}$) iff one of 
these conditions holds (in which case all of them hold).

\begin{theorem} Let $K$ be a field of characteristic $2$, $\theta$ a Tits endomorphism, $k =K^{\theta}$ and 
$L$ a $k$-subspace of $K$ with $1 \in L$ and $k[L] =K$. Let $M(U,\tau)$ be the generalized Suzuki Moufang set 
as defined in 5.5.  Then the 
following statements are equivalent:
\begin{enumerate}
\item $\theta$ is surjective.
\item $K$ is perfect.
\item $U$ has a Suzuki partition. 
\item $H$ acts transitively on $Z(U)^{\#}$.
\item For every $a \in U^{\#}$ there is a $b \in Z(U)^{\#}$ with
$a \in V_b$.
\end{enumerate}
\end{theorem}
\begin{bew}
It is clear that (a) and (b) are equivalent since $\theta^2$ is the Frobenius endomorphism.

Moreover, we have $(0,x) \tau \mu_{(a,b)} = (0, N(a,b)^{\theta} x)$ for $a,b,x \in L, (a,b) \ne (0,0)$.
Since $H=\langle \tau \mu_{(a,b)}~|~(a,b) \in U^{\#}\rangle$, one sees immediately that (d) implies (a).

Since $N(0,a) =a^{\theta}$ and thus $(0,x)\tau\mu_{(0,a)} =(0,a^2 x)$ for $a,x \in L, a\ne 0$,
we conclude that (d) follows from (b).

We show (e) implies (a): Let $a \in L^{\#}$. Then there is $b \in L$ with
$\mu_{(1,a)} =\mu_{(0,b)}$. This implies $1 +a^{\theta} +a = N(1,a) =N(0,b) = b^{\theta}$,
thus $a = (b+a+1)^{\theta} \in k$. Since $L$ generates $K$ as a ring, this implies 
$K=k$. 

Next we show that (c) follows from (a). If $a,b \in L$ with $a \ne 0, b \ne 0, b \ne a^{1+\theta}$, 
then $(a,b) =-\sim (0,s) +\sim (0,t)$ with $t =({a \over b})^{\theta^{-1}}, s = a + t$.

Finally we show that (c) implies (e). We only have to show that if $(a,b) = -\sim (0,s) +\sim (0,t)$ with 
$s,t \ne 0, s \ne t$, then there is a $c \in L$ with $\mu_{(a,b)} = \mu_{(0,c)}$. One has 
$(a,b) = (s+t, (s+t) t^{\theta})$ and thus $$N(a,b) = (s+t)^{2+\theta} + (s+t)^2 t^{\theta} + t^2 (s+t)^{\theta}
= (s^2 +t^2) s^{\theta} + s^{\theta} t^2 + t^{2+\theta}$$
$$ = s^{2+\theta} +t^{2+\theta} = (s^{1+\theta} +t^{1+\theta})^{\theta}
=N(0,s^{1+\theta} + t^{1+\theta}),$$ which yields (e)
\end{bew}

\end{document}